\documentclass[]{{birkart}}

\usepackage{graphicx}
\usepackage{amsfonts}
\usepackage{amsmath,amscd}
\usepackage{amsthm}

{ \theoremstyle{plain}
\newtheorem{theorem}{Theorem}
\newtheorem{lemma}{Lemma}
\newtheorem{corollary}{Corollary}
\newtheorem{proposition}{Proposition}
}

{ \theoremstyle{definition}
\newtheorem{definition}{{\bf Definition}}}

{ \theoremstyle{definition}

}

{ \theoremstyle{definition}
\newtheorem{remark}{Remark}
}

\begin{document}

\def\scal#1#2{\langle #1, #2\rangle}
\def\M{\mathcal{M}}
\def\N{{\gamma}}
\def\RF{\dashv}
\def\R#1{\mathbb{R}^{#1}}
\def\simm{\;\stackrel{\M}{\sim}\;}
\def\definition{\smallskip\noindent{DEFINITION.\ }}


\noindent
\textbf{On the Gauss map of embedded minimal tubes}
\footnote{This paper was supported by Sankt-Peterburg University grant,
project 95-0-1.9-34}

\bigskip
\leftline{{I.M.~Reshetnikova and V.G.~Tkachev}}

\begin{small}

\bigskip
\noindent
{\bf Abstract.} A surface is called a tube if its level-sets with respect to some coordinate function (the axis of the surface) are compact. 
Any tube of zero mean curvature has an invariant, the so-called flow vector. We study how the geometry of the Gaussian image of a higher-dimensional minimal tube $\M$ is controlled by the angle $\alpha(\M)$ between the axis and the flow vector of $\M$.
We prove that the diameter of the Gauss image of $\M$ is at least
$2\alpha(\M)$ if the angle  $\alpha(\M)$ is positive. As a consequence we derive an estimate on the length of
a two-dimensional minimal tube $\M$ in terms of $\alpha(\M)$ and the total Gaussian curvature of $\M$.


\end{small}

\section{Introduction}

Let $M$ be a Riemannian orientiable manifold of dimension
$(n-1)$ with $n\geq 3$ and $x=x(m):M\to \R{n}$ be an embedded
surface $\M$. Further we identify $\M$ with $x(M)$.
The following definition is due \cite{M1}

\begin{definition}
Let $t(\M)\subset\R{1}$ be an open interval. A surface $\M$ is
called a {\it tube} (\textit{tubular}) with the axis $Ox_n$ (or a tube in
the $e_n$-direction) with the projection interval $t(\M)$ if

1) $\forall t\in t(\M )$ the cross-sections
$\Sigma_{t}=\M\bigcap\Pi_{t}$, $\Pi_{t}=\{ x\in\R{n} : x_n=t\}$,
are nonempty compact sets lying in the interior of $\M$ (other
words, the preimage $x^{-1}(\Sigma_t)$ is a compact subset of
$M$);

2) $\forall t_1, t_2\in t(\M)$ any portion of $\M$ situated
between $\Pi_{t_1}$ and $\Pi_{t_2}$ is compact.

In this case, the length $|t(\M)|$ of the projection interval $t(\M)$ is
called the {\it life-time} of $\M$.
\end{definition}

In the present paper we study the geometry of the Gaussian image of the minimal tubes (that is tubes having zero mean curvature).
We notice that the previous definition does not impose
a priori restrictions on the topological structure of $\M$.
Important results on two-dimensional minimal tubes have been
obtained by J.C.C.Nitsche in \cite{N1}, \cite{N2}. We should
mention also the substantial papers of W.H.Meeks and
B.White~\cite{MW}, \cite{MW1} and Y.Fang \cite{YF} where
minimal tubes with convex sections $\Sigma_t$ have been studied.

The simplest example of a minimal tube in higher dimensions are the
rotationally symmetric minimal surfaces (the so-called
$(n-1)$-dimensional catenoids). In \cite{MV} (see also \cite{MT0}) V.M.Mik\-lyu\-kov and
A.D.Vedenyapin obtained finiteness of the life-time of {\it
every} minimal $(n-1)$-dimensional tube in $\R{n}$ in the case
$n\geq 4$.

The situation changes in the two-dimensional case. An
example of a tube with infinite life-time is the standard
catenoid. On the other hand, there are properly embedded singly
periodic minimal surfaces $\M$ constructed by B.~Riemann in
\cite{R} and their generalizations given in \cite{CHM} which
produce by quotient the minimal tube $\M_1=\M/{\bf Z}$ of finite
life-time. The latter means that $\M_1$ can not be a proper part
of any larger minimal tube.

In the recent paper \cite{T} the second author proposed a new approach to
the problem of determining  whether given a minimal tube be of
finite or infinite life-time. The main tool is the notion of the
flow vector of a minimal tube $\M$ which is defined below. In the
present paper we show (Theorem 1 below) that the length of the
Gaussian image of section $\Sigma_t$ can be described in terms
of slope of the flow vector to the axis of $\M$.

Let $t$ be a regular value of the coordinate function $x_n$.
Then $\Sigma_{t}$ splits in finite union of compact
$(n-2)$-dimensional connected submanifolds of $\M$. Let $\nu$ be
the unit exterior normal to $\Sigma_{t}$ as a boundary of
$\M\cap\{x_n<t\}$. In particularly, we have $\scal{\nu}{e_n}> 0$
everywhere in $\Sigma_t$.

\smallskip
\begin{definition}
A union $\Sigma=\bigcup_{i=1}^k\Sigma^i$ of components of
$\Sigma_{t}$ with the induced orientation is called a {\it
cycle}. The linear functional $F(e)$ defined by
$$
F_{\Sigma}(e)=\int\limits_{\Sigma}
\scal{e}{\nu}:\R{n} \rightarrow \R{1}
$$
generates the dual element $J(\Sigma)\in\R{n}$ such that
\begin{equation}
F_{\Sigma}(e)=\scal{e}{J(\Sigma)}\, .
\label{e1}
\end{equation}
\end{definition}

Then
\begin{equation}
J_n(\Sigma)\equiv\scal{J(\Sigma)}{e_n}=\int\limits_{\Sigma}
\scal{e_n}{\nu}= \int\limits_{\Sigma} |e_n^\top|>0,
\label{posit}
\end{equation}
and it follows that $J(\Sigma)\ne 0$. Here and subsequently we
use the notation $a^V$ for the orthogonal projection of
$a\in\R{n}$ onto a subspace $V\subset\R{n}$ and by $T=T_m\M$ we
denote the tangent space of the surface $\M$ at $m\in \M$.

\smallskip
\definition
Two oriented cycles $\Sigma'$ and $\Sigma''$ are {\it
equivalent} in $\M$, or $\Sigma'\simm\Sigma''$, if there exists
an open subset $D\subset \M$ such that $\partial D
=(-\Sigma')\bigcup\Sigma''$ (here $-\Sigma'$ is the opposite
oriented cycle to $\Sigma$). This notion is actually the
oriential bordism equivalence (see \cite[\S~7]{Hirsh}).  A
connected cycle $\Sigma$ is called {\it simple} if it is
equivalent to zero cycle in the hyperplane $\Pi_t$. Other words,
it is a boundary of an open subset of $\Pi_t$.

\begin{proposition}
We have for  $\Sigma=\bigcup_{i=1}^k\Sigma^i$
$$
J(\Sigma)=J(\Sigma^1)+\ldots+J(\Sigma^k)\, .
$$
Moreover, if $\Sigma'\simm\Sigma''$ then
$J(\Sigma')=J(\Sigma'')$.
\label{pr1}
\end{proposition}

\begin{proof}
The first property is direct consequence of the above
definitions. To prove the second one we recall that
 all coordinate functions of minimal immersion are
harmonic~\cite{KN}. Let $D\subset \M$ be an open set in the Definition~2
such that $\partial D =(-\Sigma')\bigcup\Sigma''$. Then for
arbitrary coordinate vector $e_k\in \R{n}$
$$
\scal{J(\Sigma'')}{e_k} -
\scal{J(\Sigma')}{e_k}=\int\limits_{\partial D} \scal{e_k}{\nu}=
\int\limits_{\partial D} \scal{\nabla f_k}{\nu}
=\int\limits_D \Delta f_k = 0\, ,
$$
where $\nabla f_k=e_k^\top$ is the gradient of $f_k=\scal{e_k}{x(m)}$.
\end{proof}

\smallskip
\definition
We call $J(\M)=J(\Sigma_t)$ to be the {\it flow vector} of the
tube $\M$.

\remark
It follows from Proposition~1 that the flow vector
$J(\Sigma_t)$ does not depend on a choice of $t\in t(\M)$. One
easy to see also that both the angle $\alpha(\M)$ between
$J(\M)$ and $e_n$ and the norm $\|J(\M)\|$ are invariants under
the action of the orthogonal subgroup of $\R{n}$ preserving the
axis $Ox_n$.

Moreover, we emphasize that the flow vector of $\M$ is a local
characteristic of $\M$ in the sense that it can be
computed if we consider only a portion of $\M$
situated between $\Pi_{t_1}$ and $\Pi_{t_1}$ for $t_1$ and $t_2$
arbitrarily close.

Let $S^{n-1}$ be the unit sphere in the Euclidean space $\R{n}$
and $d(E)$  be the spherical diameter of a set $E\subset S^{n-1}$.
By $\gamma:\M\to S^{n-1}$ we denote the Gaussian map of $\M$,
where $\gamma(m)$ is the unit normal at $m\in \M$; by $\gamma(E)$
we denote the Gaussian image of a set $E\subset \M$.

Our main result is the following lower estimate of the diameter of
the Gaussian image.

\begin{theorem}
Let $\M$ be an embedded minimal tube in $\R{n}$;
$\Sigma\subset\Sigma_t$ be a simple cycle with the flow vector
$J(\Sigma)$. Then the diameter of $\gamma(\Sigma)$ satisfies
\begin{equation}
d(\gamma(\Sigma))\geq 2\alpha(\Sigma)\, ,
\label{e2}
\end{equation}
where $\alpha(\Sigma)$ is the angle between $J(\Sigma)$ and $e_n$.
\label{th1}
\end{theorem}

As a consequence in Section~\ref{sec4} we obtain the upper
estimate on the life-time of minimal tubes of finite total
Gaussian curvature.

\begin{theorem}
Let $\M$ be a two-dimensional minimal tube in $\R{3}$ of finite
total Gaussian curvature $-G(\M)$. If $\alpha(\M)>0$  then
$\M$ has finite life-time and
\begin{equation}
|t(\M)|\leq \|J(\M)\|\;G(\M)\frac{\cos\alpha(\M)}{16\alpha^2(\M)}.
\label{tube}
\end{equation}
\label{th2}
\end{theorem}

\begin{corollary}
Let $\M$ be a two-dimensional minimal tube in $\R{3}$ with
univalent Gaussian map. Then $\M$ has finite life-time provided
that $\alpha(\M)>0$.
\label{cor1}
\end{corollary}

Now we indicate the main idea of the proof of Theorem~\ref{th2}.
In this case $\dim \Sigma=1$ and it follows that all
one-dimensional cycles are simple (see \cite[\S~7]{Hirsh}).
Moreover, (\ref{e2}) implies that the Gaussian image of every
section $\Sigma_t$ is uniformally `large' provided the angle between
$J(\M)$ and $e_n$ is strictly positive. On the other hand, in
the two-dimensional case $\dim \M=2$ the Gaussian map is
conformal and (\ref{e2}) yields that $\M$ must be a surface of
hyperbolic conformal type. The final step is to use the
connection between the conformal module of minimal tube and
its life-time value.

We notice that Theorem~\ref{th2} fails if we drop the requirements
of finiteness of the total Gaussian curvature. Really, in the
previous paper \cite{T1} we have constructed the corresponding
examples by using the suitable Weierstrass representation for
minimal tubes. Namely, given arbitrary $\alpha(\M)>0$ there exists
a properly embedded minimal tube of infinite life-time.

The author thanks Yi Fang for helpful discussions concerning the
subject of the paper.

\section{Preliminary facts}

By $\Lambda(\R{n})$ and $\Lambda^k(\R{n})$ we denote the exterior
algebra of $\R{n}$ and the subspace of all $k$-form respectively.
We specify an orthonormal basis $\{e_k\}^n_{k=1}$ in $\R{n}$ and
by $\omega=e_1\wedge\ldots\wedge e_n$ we denote the volume-form
of $\R{n}$. We write $a\simeq b$ if $a=\pm b$.

Given $u\in\Lambda(\R{n})$ we define the inner product $u\RF
\cdot$ on $\Lambda(\R{n})$ by
\begin{equation}
\scal{x}{u\RF y}\equiv\scal{u\wedge x}{y}\, , \quad \forall x,
y\in\Lambda(\R{n}) \,.
\label{e3}
\end{equation}
Then the Hodge $*$-perator:
$\Lambda^k(\R{n})\rightarrow\Lambda^{n-k}(\R{n})$,
can be written for every $k$-form $x$ by
\begin{equation}
*x=x\RF \omega\, .
\label{e4}
\end{equation}

The following facts are elementary and can be found in \cite{St}.

(i) $**x=(-1)^{k(n-k)}x\simeq x\, ,\quad \forall x\in
\Lambda^k(\R{n})$;

(ii) $x\wedge *y=\scal{x}{y}\omega\, ,\quad
\forall x\in\Lambda^k(\R{n}), \forall
y\in\Lambda^{n-k}(\R{n})\,;$

(iii) $\scal{*x}{*y}=\scal{x}{y}\, ,\quad\forall x,
y\in\Lambda^k(\R{n})$. \\

Let $V\subset\R{n}$ be an oriented $k$-dimensional subspace
and $v_1,\ldots,v_k$ be an orthonormal basis of $V$. By
$$
\sigma(V)\equiv v_1\wedge\ldots\wedge v_k\,
$$
we denote the volume form of $V$. Further we use the operator
$$
\pi_V(\xi)=*\left(\sigma(V)\wedge\xi\right)\, .
$$

Let $V\subset\R{n}$ be an oriented hyperspace, ${\rm
dim} V=n-1$. If $\xi$ is a unit vector orthogonal
to $V$ (i.e. $\xi\in V^\bot$)  then it follows from~(ii)
\begin{equation}
\xi\simeq *\sigma(V).
\label{e5}
\end{equation}

\begin{lemma}
For all $a\in \Lambda^r(\R{n})$, $b\in \Lambda^k(\R{n})$
one holds
$$
a\RF(*b)=*(b\wedge a)\, .
$$
\label{lem1}
\end{lemma}

\begin{proof}
Let $\xi\in \Lambda^{n-k-r}(\R{n})$ be chosen arbitrary. Then
$$
\scal{a\RF (*b)}{\xi}=\scal{*b}{a\wedge
\xi}=\scal{b\RF \omega}{a\wedge \xi}=
$$
$$
=\scal{\omega}{b\wedge a\wedge \xi}=\scal{(b\wedge a)\RF \omega}
{\xi}=\scal{*(b\wedge a)}{\xi}\, .
$$
Then by duality we have $*(b\wedge a)=a\RF (*b)$ and the lemma is
proved.

\end{proof}

\begin{lemma}
Let $V\subset\R{n}$ be a subspace,  $\dim  V=k$, and
$V^\bot$  be its orthogonal complement. Then for all
$x\in\R{n}$
\begin{equation}
x^{V^\bot}=\sigma(V)\RF \left(\sigma(V)\wedge x\right)\, .
\label{e6}
\end{equation}
\label{lem2}
\end{lemma}

\begin{proof}
We choose $v_1,\ldots, v_k$ to be an orthonormal basis of $V$ and
consider its complement in $\R{n}$: $w_1,\ldots,w_{n-k}$. Given
arbitrary $x\in\R{n}$ we have
$x=x_1 v_1+\ldots+x_k v_k+y_1w_1+\ldots+y_{n-k}w_{n-k}$.
Then for every $v_i$
$$
\scal{\sigma(v)\RF \left(\sigma(v)\wedge x\right)}{v_i}=
\scal{\sigma(v)\wedge x}{\sigma(v)\wedge v_i}=0\, .
$$
On the other hand,
$$
\scal{\sigma(v)\RF \left(\sigma(v)\wedge x\right)}{w_\alpha}=
\scal{\sigma(v)\wedge x}{\sigma(v)\wedge w_\alpha}=
$$
$$
=\sum\limits_{j=1}^{n-k} \scal{\sigma(v)\wedge w_j}{\sigma(v)
\wedge w_\alpha}=\sum\limits_{j=1}^{n-k} y_j \delta_{j\alpha}=
y_\alpha\, .
$$

Hence, by the definition we have the identity
$$
\sigma(v)\RF \left(\sigma(v)\wedge x\right)=
y_1w_1+\ldots+y_{n-k}w_{n-k}=x^{V^\bot}\, ,
$$
which proves the lemma.

\end{proof}

\section{Proof of Theorem 1}

In this section by $\Pi=\Pi_0$ we denote a hyperspace $x_n=0$ in
$\R{n}$ and by $T\equiv T_m\Sigma$ --- the tangent space to the
section $\Sigma$ being considered as submanifold of $\Pi$. Let
$\N=\N(m)$ be the unit normal to $\M$ at $m$. We specify an
orthonormal basis $\tau_1,\ldots,\tau_{n-2}$ in $T$ and by
$\tau\equiv\sigma(T)$ denote the volume form of $T_m\Sigma$.

We need the following auxilliary assertion

\begin{lemma}
Let $\xi$ and $\eta$ be two unit vectors such that
$\xi,\eta,\tau_1,\ldots,\tau_{n-2}$ form oriented orthonormal
basis of $\R{n}$. Then for every $q\in\R{n}$
$$
\scal{q}{\xi}\simeq\scal{\eta}{\pi_T(q)}.
$$
\label{lem4}
\end{lemma}

\begin{proof}
We have from (\ref{e5}) that
$$
\xi\simeq*(\eta\wedge \tau),
$$
and by virtue of (iii) and Lemma~\ref{lem1} we obtain
$$
\scal{q}{\xi}\simeq\scal{q}{*(\eta\wedge \tau)}
\simeq\scal{q}{*(\tau\wedge \eta)}\simeq\scal{q}{\eta\RF*\tau}\simeq
$$
$$
\simeq\scal{\eta\wedge q}{*\tau}
\simeq\scal{\eta}{q\RF*\tau}
\simeq\scal{\eta}{*(\tau\wedge q)}\simeq\scal{\xi}{\pi_T(q)}
$$
and the lemma is proved.

\end{proof}

By the Sard's theorem and regularity of $t$ we conclude that
$\Sigma$ is a smooth submanifold of $\Pi_t$. Assume that
$\eta=\eta(m)$ is the unit normal vector field to $\Sigma$ in
$\Pi_t$ oriented such that the pair $(T;\eta)$ is an oriented
basis of $\Pi$.  Then by Lemma~\ref{lem2} we have for
every $q\in\R{n}$
$$
\int\limits_{\Sigma}\scal{\pi_T(q)}{e_n}=
\int\limits_{\Sigma}\scal{*(\tau\wedge q)}{e_n}\simeq
\int\limits_{\Sigma}\scal{\tau\wedge q}{*e_n}\simeq
$$
\begin{equation}
\simeq\int\limits_{\Sigma}\scal{\tau\wedge q}{\tau\wedge\eta}
\simeq\int\limits_{\Sigma}\scal{q}{\tau\RF(\tau\wedge\eta)}
\simeq\int\limits_{\Sigma}\scal{q}{\eta}
\label{eq1}
\end{equation}

To show that in fact the last integral  vanishes we observe
that by the definition, the simple cycle $\Sigma$ is the boundary
of some open subset $\Omega\subset \Pi_t$. The by the Stokes'
formula we obtain
\begin{equation}
\int\limits_{\partial\Omega}\scal{q}{\eta}=
\int\limits_{\Omega} {\rm div}\;q=0.
\label{need}
\end{equation}

Thus (\ref{eq1}) yields the following identity
\begin{equation}
\int\limits_{\Sigma}\scal{\pi_T(q)}{e_n}=0.
\label{eq2}
\end{equation}

Choose $q\ne0$ arbitrarily such that the equality
$\scal{q}{J(\Sigma)}=0$ holds. Then taking into account the
mutual orthogonality of $\N$, $\nu$ and the tangent space
$T_m\Sigma$, we obtain from (\ref{e1}) and Lemma~\ref{lem4}
$$
0=\int\limits_{\Sigma}\scal{q}{\nu}=
\int\limits_{\Sigma}\scal{\pi_T(q)}{\N}.
$$
Hence, we conclude from (\ref{eq2}) that
\begin{equation}
\int\limits_{\Sigma}\scal{\pi_T(q)}{\N\pm e_n}=0.
\label{eq3}
\end{equation}

By virtue of regularity of $t$, the expressions $\N\pm e_n$ does
not vanish everywhere in $\Sigma$. It follows that along $\Sigma$
the vector fields
$$
v_{\pm}(m)=\frac{e_n\pm \N(m)}{\|e_n\pm \N(m)\|}
$$
are well-defined.

Using the mean value theorem  we deduce from  (\ref{eq3}) that
there exist two points $m_-$ and $m_+$ in $\Sigma$ such that
\begin{equation}
\scal{\pi_{T_{\pm}}(q)}{v_\pm}=0,
\label{eq4}
\end{equation}
where $T_{\pm}=T_{m_\pm}\M$ ¨ $v_\pm=v_\pm(m_\pm)$.

Now we observe that the set $(v_-,v_+,\tau_1,\ldots,\tau_{n-2})$
forms an orthonormal basis $\R{n}$. Thus, applying
Lemma~\ref{lem4} to  (\ref{eq4}) we obtain
$$
0=\scal{\pi_{T_{+}}(q)}{v_+}\simeq
\scal{v_-}{q}
$$
at $m_-$ and similarily at $m_+$:
$$
0=\scal{\pi_{T_{-}}(q)}{v_-}\simeq
\scal{v_+}{q}.
$$

By getting rid of the denominator in the definition of the
vectors $v_\pm$ we arrive at
$$
\scal{q}{e_n\pm \N(m_\pm)}=0.
$$
It follows that
$$
\scal{\N(m_+)-\N(m_-)}{q}=2q_n,
$$
where $q_n=\scal{q}{e_n}$. Finally,
applying the Cauchy`s integral inequality in the last identity
yields
\begin{equation}
\|\N(m_+)-\N(m_-)\|\geq \frac{2q_n}{\|q\|}.
\label{eq5}
\end{equation}

Taking into account that $\N(m_\pm)$ are points on the unit
sphere lying in the Gaussian image of $\Sigma$,  we obtain
$\|\N(m_+)-\N(m_-)\|=2\sin\frac{\beta}{2}$, where $\beta$
is the angle between $\N(m_+)$ and $\N(m_-)$. It follows
by the definition of the spherical diameter of $\N(\Sigma)\subset
S^{n-1}$ that $d(\N(\Sigma))\geq\beta$. Thus, by virtue of
(\ref{eq5}) we conclude that
$$
d(\N(\Sigma))\geq 2 \arcsin \frac{q_n}{\|q\|}.
$$

To find the maximum of the right part of the last expression we
assume $\alpha(\Sigma)$ to be equal the angle between
$J(\Sigma)$ and $e_n$. Then, by orthogonality of $q$ to the
flow vector $J(\Sigma)$ one easily sees that
$$
\max_{q\bot J(\Sigma)}\frac{q_n}{\|q\|}=
\sin\alpha(\Sigma).
$$
Hence,
$$
d(\N(\Sigma))\geq 2 \alpha(\Sigma),
$$
and Theorem~\ref{th1} is proved completely.

\remark
We notice that in the two-dimensional case the assertion
of Theorem~\ref{th1} is still true provided $\M$ is properly immersed
minimal tube. To check this fact we observe that the unique
place in the proof of the theorem where we essentially needed the
embeddedness hypothesis is formula (\ref{need}). The validness of this
formula in the two-dimensional immersed case is a direct
consequence of the Green integration formula.


\section{Applications to two-dimensional minimal tubes}
\label{sec4}

To prove Theorem~\ref{th2} we need some terminology from
potential theory.

Let us consider an embedded minimal hypersurface $\M$ in $\R{n}$
which is a tube in $e_n$-direction. Given $t_1$, $t_2$ from the
interval $t(\M)$ we notice by $\M(t_1;t_2)$ the portion of $\M$
situated in the slab $t_1<x_n<t_2$. Then the quantity
$$
{\rm cap}\;\M(t_1;t_2)=\inf \int\!\!\int_{\M(t_1;t_2)}|\nabla \varphi|^2,
$$
where the infimum is taken over all Lipschitzian functions
$\varphi(m)$ on $\M(t_1;t_2)$ such that $\varphi(m)=0$ on
$\Sigma_{t_1}$ and $\varphi(m)=1$ on $\Sigma_{t_2}$ is called
the {\it capacity} of $\M(t_1;t_2)$.

Let $\Gamma$ be a family of locally rectifiable curves
$\gamma\subset \M$ and $\rho(m)\geq 0$ be a Baire function with
the property
$$
\int_\gamma \varphi(z)\,ds\geq 1,
$$
for every  $\gamma\in \Gamma$. The infimum
$$
{\rm mod} \;\Gamma=\inf \int\!\!\int_{\M(t_1;t_2)} \rho^2(m)
$$
over all such $\rho(m)$ is called the {\it module }
of the family $\Gamma$.

If $\dim \M=2$ the following connection between the capacity of
$\M(t_1;t_2)$ and the module of the family $\Gamma(t_1;t_2)$ of
all curves which connect two boundary components of
$\M(t_1;t_2)$ is well-known
$$
{\rm mod} \;\Gamma(t_1;t_2)={\rm cap}\;\M(t_1;t_2)
$$
(see for the Euclidean case \cite{Fug} and for the Riemannian
case \cite{Mikl} respectively).

In his paper \cite{M1} V.M.Miklyukov has studied the
higher-dimensional minimal tubes in $\R{n}$ and has established the
following connection between the capacity of $\M(t_1;t_2)$ and
its life-time
\begin{equation}
{\rm cap}\;\M(t_1;t_2)=\frac{t_2-t_1}{\scal{J(\M)}{e_n}}.
\label{mik}
\end{equation}

{STEP 1}. \
Let us first assume that $\M$ be a two-dimensional
embedded minimal tube to be homeomorphic to an annulus. Then from
(\ref{mik}) we obtain
\begin{equation}
{\rm mod}\;\Gamma(t_1;t_2)=\frac{t_2-t_1}{\scal{J(\M)}{e_3}}.
\label{rik}
\end{equation}

Since the Gaussian map of a minimal
surface is conformal, we have for every tangent vector $X\in T_mM$
$$
\|d\N_m (X)\|=\lambda(m) \|X\|,
$$
and the Gaussian curvature is $K(m)=-\lambda^2(m)$.

Given an arbitrary $t\in t(\M)$ we write by change coordinates
formula
$$
\int\limits_{\Sigma_t} \lambda\, ds\geq
\int\limits_{\N(\Sigma_t)} ds_1\geq 2 d(\N(\Sigma_t)),
$$
where $ds_1$ is the metric element on the unit sphere.
By virtue of Theorem~\ref{th1} we obtain
\begin{equation}
\int\limits_{\Sigma_t} \lambda\, ds\geq
4\alpha(\M),
\label{eq6}
\end{equation}
where $\alpha(\M)$ is the angle between the flow vector $J(\M)$
and $e_3$.

Substituting $\rho(m)=\lambda(m)$ in the definition of the module,
we get from (\ref{eq6})
$$
{\rm
mod}\;\Gamma(t_1;t_2)\leq\frac{1}{16\alpha^2(\M)}
\int\!\!\int_{\M(t_1;t_2)}(-K),
$$
and by (\ref{rik}) we arrive at
$$
t_2-t_1\leq\frac{J_3(\M)G(\M(t_1;t_2))}{16\alpha^2(\M)}.
$$
Then the required estimate in the annulus case yields from the
arbitrariness of $t_1$ and $t_2$.

In order to prove the general case we need the following elementary fact.

\begin{lemma}
Let $v_1$, \ldots $v_l$ be a system of nonzero vectors from
$\R{n}$ such that for some $e\in\R{n}$ one holds
$\alpha_i\leq\pi/2$, where $\alpha_i$ is the angle between $v_i$
and $e$. Let $v=\sum_{i=1}^{l}v_i$.  Then we have for the angle
$\alpha$ between $v$ and $e$:
$$
\alpha\leq\max\{\alpha_1,\ldots,\alpha_l\}.
$$
\label{convex}
\end{lemma}

\begin{proof}
By virtue of $\scal{v_i}{e}\geq 0$ we notice that
$\alpha_i\leq\pi/2$ and by the triangle inequality obtain
$$
\cos\alpha=\frac{\scal{v}{e}}{\|v\|}\geq
\frac{\scal{v}{e}}{\sum_{i=1}^{l}\|v_i\|}=
\frac{\sum_{i=1}^{l}\scal{v_i}{e}}{\sum_{i=1}^{l}\|v_i\|}=
$$
$$
=\frac{\sum_{i=1}^{l}\|v_i\|\cos\alpha_i}{\sum_{i=1}^{l}\|v_i\|}\geq
\min_{1\leq j\leq l}\cos\alpha_i=\cos(\max_{1\leq j\leq l}\alpha_j),
$$
as required.

\end{proof}

STEP 2.\
Let $\M$ be a properly embedded minimal tube in $\R{3}$ of
general topological structure.
First we notice that for
arbitrary closed subinterval $[\alpha;\beta]\subset t(\M)$ there
exist at most finitely many points $m\in\M(\alpha;\beta)$ such
that $\gamma(m)=\pm e_3$. Really, the coordinate function
$f_3(m)=\scal{e_3}{x(m)}$ is harmonic on $\M$ and it follows
that the critical set $H=\{m\in\M:\nabla f_3\equiv e_3^\top=0\}$
has no accumulation points inside of $\M$. Other words, for any
compact part $\M(\alpha;\beta)$ the set $H$ is finite.

Let $c_1<c_2<\ldots c_{k-1}$ are all the values of $f_3(m)$ when
$m$ runs $H$ and $\alpha=c_0$, $\beta=c_k$. Then by the Morse
theory, every part of $\M_i=\M(c_{i-1},c_i)$, $i=1,\ldots, k$,
is a union of annuli $\M_i^1$, \ldots $\M_i^{\ell_i}$.

From positiveness  of the third coordinate of the
flow vector (\ref{posit}) and Proposition~1 we have
$J_3(\M_i^j)\leq J_3(\M)$. Applying Step~1 we obtain
\begin{equation}
c_i-c_{i-1}\leq \frac{J_3(\M)G(\M_i^{j})}{16\alpha^2(\M_i^j)},
\qquad 1\leq j\leq \ell_i.
\label{eer}
\end{equation}

On the other hand, by virtue of Proposition~1
$$
J(\M)=J(\M_i)=J(\M_i^1)+\ldots+J(\M_i^{\ell_i}).
$$
Let the index $j_0$ corresponds to the maximum angle
$\alpha(\M_i^j)$ when $i$ is fixed. Then applying
Lemma~\ref{convex} to $e=e_3$ and $v_j=J(\M_i^j)$ we obtain
$$
\alpha(\M_i^{j_0})\geq \alpha(\M_i)=\alpha(\M).
$$
Hence, by (\ref{eer}) and positiveness of the absolute total
Gaussian curvature $G$
$$
c_i-c_{i-1}\leq
\frac{J_3(\M)G(\M_i^{j_0})}{16\alpha^2(\M)}\leq
\frac{J_3(\M)G(\M_i)}{16\alpha^2(\M)}.
$$
Summing of the last inequalities over all $i=1,\ldots,k$ we
arrive at
$$
\beta-\alpha\leq
\frac{J_3(\M)G(\M(c_0;c_k))}{16\alpha^2(\M)}\leq
\frac{J_3(\M)G(\M)}{16\alpha^2(\M)}=
\|J(\M)\|G(\M)\frac{\cos\alpha(\M)}{16\alpha^2(\M)}.
$$
By arbitrariness of subinterval $[\alpha;\beta]\subset t(\M)$
we obtain the assertion of Theorem~2.


\bigskip
\bigskip
\noindent
Department of Mathematics

\noindent
Volgograd State University

\end{document}